%% file: main.tex
\documentclass[12pt]{article}
\usepackage[reqno]{amsmath}
\usepackage{authblk}
\usepackage{amssymb, enumerate,amsthm}
\usepackage{graphicx}
\usepackage{tikz}
\usetikzlibrary{calc}
\usetikzlibrary{decorations.markings}
\usetikzlibrary{decorations.text}
\usepackage{array}
\usepackage{float}
\usepackage{wrapfig}
\usepackage{fullpage}
\usepackage{url}
\usepackage{hyperref}
\usepackage{nicefrac}
\usepackage{blkarray}
\usepackage{dcolumn,booktabs}
\newcolumntype{d}[1]{D{.}{.}{#1}}
\newcommand\mc[1]{\multicolumn{1}{c}{#1}}
\usepackage{nicematrix}
\usepackage{enumitem}
\usepackage[capitalise]{cleveref}
\usepackage{color,colortbl}

\newcommand{\doi}[1]{\href{https://doi.org/#1}{\texttt{doi:#1}}}

\expandafter\ifx\csname pdfpagebox\endcsname\relax\else
\fi
 
\newtheorem{theorem}{Theorem}[section]
\newtheorem{lemma}[theorem]{Lemma}

\newcommand\GM[1]{\mathrm{GM}(#1)}
\newcommand\KS[1]{\mathrm{KS}(#1)}

\definecolor{gcolour}{RGB}{230,97,0}
\definecolor{kcolour}{RGB}{93,58,155}

\newcommand\email[1]{\texttt{#1}}

\title{Cubic graphs with no eigenvalues in the interval $(-1,1)$ }
\author[1]{Krystal Guo}
\author[2]{Gordon F. Royle}
\affil[1]{Korteweg-de Vries Institute for Mathematics, University of Amsterdam. (\email{k.guo@uva.nl}).}
\affil[2]{Department of Mathematics and Statistics, The University of Western Australia. (\email{gordon.royle@uwa.edu.au}).}
\begin{document}
\maketitle

\begin{abstract}
We give a complete characterisation of the cubic graphs with no eigenvalues in the open interval $(-1,1)$. There are two infinite families, one due to Guo and Mohar [Linear Algebra Appl. 449:68--75] the other due to Koll\'ar and Sarnak [Communications of the AMS. 1,1--38], and $14$ ``sporadic'' graphs on at most $32$ vertices. This allows us to show that $(-1,1)$ is a maximal spectral gap set for cubic graphs. Our techniques including examination of various substructure and an application of the classification of generalized line graphs.

\vspace{5pt}
\noindent \emph{MSC2020:} Primary 05C50; Secondary 05C76.

\noindent \emph{keywords: graph eigenvalues, graph classification, graph spectra, spectral gap set}
\end{abstract}


\input intro2

\input examples

\input proof

\input{conclusion}

\bibliographystyle{plain}
\bibliography{median}

\end{document}

%% file: intro2.tex
\section{Introduction}


We show that with the exception of  two infinite families and 14 sporadic examples, all connected cubic graphs have an eigenvalue in the open interval $(-1,1)$.

This work was motivated by results from two related topics in spectral graph theory, namely the study of \emph{spectral gap sets} for cubic graphs, and the study of \emph{median eigenvalues} for chemical graphs.

If $G$ is a cubic graph then its eigenvalues (that is, the eigenvalues of its adjacency matrix) lie in the interval $[-3,3]$. Koll\'ar and Sarnak \cite{KolSar2021} (but see also \cite{kollar1}) define a \emph{spectral gap set}, or just \emph{gap set} to be an open subset $\mathcal I \subset [-3,3]$ such that there are infinitely many cubic graphs with no eigenvalues in $\mathcal I$. A gap set $\mathcal{I}$ is \emph{maximal} if there is no gap set $\mathcal{I}'$ such that $\mathcal{I} \subset \mathcal{I}'$. 

Koll\'ar and Sarnak observe that several important results in graph theory are equivalent to identifying maximal gap intervals. For example, the Alon-Boppana theorem (see \cite{Alo1986, HooLinWig2006, Cio2006,Moh2010}) and the existence of cubic Ramanujan graphs  \cite{Chi1992} show that $(2\sqrt{2},3)$ is a maximal gap interval.

Among several other results, Koll\'ar and Sarnak \cite{KolSar2021} exhibit an infinite family of cubic planar graphs with no eigenvalues in $(-1,1)$ thereby showing that $(-1,1)$ is a gap interval (in fact they show that it is a maximal gap interval). The fact that $(-1,1)$ is a gap interval had previously been shown by Guo and Mohar \cite{GuoMoh2014} who found an infinite family of  (non-planar) cubic graphs with no eigenvalues in $(-1,1)$.

We show that these two infinite families, along with $14$ small graphs on at most $32$ vertices, form the complete list of cubic graphs with no eigenvalues in $(-1,1)$. As a result of this characterization, we can then show that $(-1,1)$ is a maximal gap set (not just a maximal gap \emph{interval}). 

A \emph{chemical graph} or a \emph{subcubic graph} is a connected simple graph with maximum degree $3$. In mathematical chemistry the eigenvalues of such a graph correspond to orbital energies in the H\"uckel model (see \cite{FowPis2010}). If we denote the eigenvalues of a chemical graph by
\[
\lambda_1 \geqslant \lambda_2 \geqslant \cdots \geqslant \lambda_n
\]
then (for even $n$) the \emph{median eigenvalues} 
\[
\lambda_{H} = \lambda_{\nicefrac{n}{2}}, \quad \lambda_{L}= \lambda_{\nicefrac{n}{2} + 1} 
\]
are the most important. They correspond to the highest occupied molecular orbital (HOMO) and the lowest unoccupied molecular orbital (LUMO) respectively, and there are several papers considering the HOMO-LUMO \emph{gap} $\lambda_H - \lambda_L$ or the HOMO-LUMO \emph{index} $R(G) = \max \left\{ |\lambda_H|, |\lambda_L| \right\}$. (When $n$ is odd, then $\lambda_{H} = \lambda_{L} = \lambda_{\nicefrac{(n+1)}{2}}$  and the first question is irrelevant.)

Fowler and Pisanski \cite{FowPis2010} consider plotting the pairs $(\lambda_H, \lambda_L)$ on the $xy$-plane, and conjecture that only a finite number of subcubic graphs lie outside the \textsl{chemical triangle} with vertices $(-1,-1)$,  $(1,-1)$ and  $(1,1)$. In fact, the Heawood graph, with $(\lambda_H, \lambda_L) = (\sqrt{2}, -\sqrt{2})$ is the only known subcubic graph not lying in the chemical triangle. Bipartite graphs have symmetric spectra and so $(\lambda_H, \lambda_L) = (-\alpha, \alpha)$ for some $\alpha$. Mohar \cite{Moh2016} exploits this to show that the Heawood graph is the only \emph{bipartite} subcubic graph with $R(G) > 1$. If we restrict further to bipartite \emph{cubic} graphs, then our results provide an alternative and more elementary proof of this result. For non-bipartite subcubic graphs the question is still open.

The remainder of the paper is structured as follows. Taking the bipartite double of a graph preserves the property of having no eigenvalues in $(-1,1)$ and so we start by determining the bipartite cubic graphs with no eigenvalues in $(-1,1)$.  First, \cref{bip4} analyses the local structure of the graph to show that the Guo-Mohar graphs are the only bipartite cubic graphs of girth $4$ with no eigenvalues in $(-1,1)$. Then, \cref{bip6} shows that any bipartite cubic graph of girth at least $6$ with no eigenvalues in $(-1,1)$ can be associated with another graph with least eigenvalue $-2$, and the famous classification of graphs with least eigenvalue $-2$ (see \cite{BusCveSei1976},\cite{CamGoeSeiShu1976}) is used to determine all the possibilities. Finally, \cref{nonbip} determines the non-bipartite cubic graphs with no eigenvalues in $(-1,1)$ by ``reversing'' the bipartite doubling process in all possible ways on each of the newly-determined cubic bipartite graphs with no eigenvalues in $(-1,1)$.





%% file: examples.tex
\newcommand{\ks}{Koll\'ar-Sarnak}
\newcommand{\gm}{Guo-Mohar}

\section{The graphs}\label{sec:graphs}

In this section we describe in detail the cubic graphs with no eigenvalues in $(-1,1)$. These fall into two infinite families and 14 ``sporadic'' graphs.

\subsection{Two infinite families}\label{families}

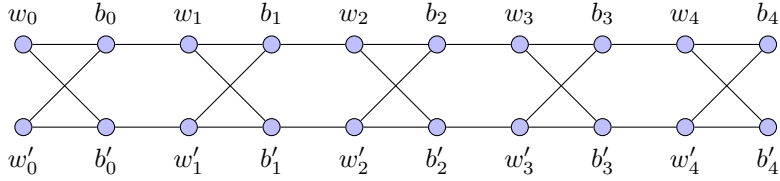
\begin{figure}
\begin{center}
\begin{tikzpicture}[scale=1.1]
\tikzstyle{vertex}=[fill=blue!25, circle, draw=black, inner sep=0.8mm]

\foreach \x in {0,1,2,3,4,5,6,7,8,9} {
\ifthenelse{\isodd{\x}}{ 
\pgfmathsetmacro\result{(\x -1)/2}
\node [vertex, label=below:{\footnotesize $b_{\pgfmathprintnumber{\result}}'$}] (v\x) at (\x,0) {};
\node [vertex, label=above:{\footnotesize $b_{\pgfmathprintnumber{\result}}$}] (w\x) at (\x,1) {}; }
{
\pgfmathsetmacro\result{\x /2 }
\node [vertex, label=below:{\footnotesize $w_{\pgfmathprintnumber{\result}}'$}] (v\x) at (\x,0) {};
\node [vertex, label=above:{\footnotesize $w_{\pgfmathprintnumber{\result}}$}] (w\x) at (\x,1) {};
}
}


\draw (v0)--(v1)--(v2)--(v3)--(v4)--(v5)--(v6)--(v7)--(v8)--(v9);
\draw (w0)--(w1)--(w2)--(w3)--(w4)--(w5)--(w6)--(w7)--(w8)--(w9);
\draw (v0)--(w1);
\draw (v1)--(w0);
\draw (v2)--(w3);
\draw (v3)--(w2);
\draw (v4)--(w5);
\draw (v5)--(w4);
\draw (v6)--(w7);
\draw (v7)--(w6);
\draw (v8)--(w9);
\draw (v9)--(w8);
\end{tikzpicture}
\caption{The base graph $B(5)$.}
\label{basefig}
\end{center}
\end{figure}

There are two infinite families, each of which exists only when the number of vertices is a multiple of $4$. There is a unique graph on $4k$ vertices in each of the families. These graphs are obtained from a \emph{base graph} which we denote $B(k)$ consisting of  $k > 1$ vertex-disjoint $4$-cycles connected in a linear fashion to from a path of $4$-cycles (as depicted for $k=5$ in Figure~\ref{basefig}).

The base graph has four vertices of degree $2$ (two at each end) and can be completed to a cubic graph in two non-isomorphic ways. The \emph{\ks} graphs are obtained by adding two edges, each joining the two vertices of degree $2$ at the same end, while the \emph{\gm} graphs are obtained by adding two edges, each joining a vertex of degree $2$ at one end to a vertex of degree $2$ at the other end. (There are two ways to do this, but the resulting graphs are isomorphic.)

More explicitly, the base graph $B(k)$ on $4k$ vertices has vertex set 
\[
 V = \{w_0, \ldots, w_{k-1}\} \cup \{w_0', \ldots, w_{k-1}'\} \cup \{b_0, \ldots, b_{k-1}\} \cup \{b_0', \ldots, b_{k-1}'\} 
\]
and edges 
\[
E = \{w_i b_i, w_i'b_i', w_ib_i', w_i'b_i \mid  i = 0,\ldots, k-1 \} \cup 
\{ b_i w_{i+1}, b_i' w_{i+1}' \mid  i = 0,\ldots, k -2 \}. 
\]
The {\ks} graph, $\KS{k}$, on $4k$ vertices is obtained from $B(k)$ by adding edges $w_0w_0'$ and $b_{k-1}b_{k-1}'$. The {\gm} graph, $\GM{k}$, on $4k$ vertices is obtained from $B(k)$ by adding edges $w_0b_{k-1}$ and $w_0'b_{k-1}'$. We see that the base graph and the {\gm} graphs are bipartite, where the bipartite colour classes are $\{w_i,w_i'\}_{i=0}^{k-1}$ and $\{b_i,b_i'\}_{i=0}^{k-1}$, while the \ks{} graphs each have four triangles. We note that $\GM{2}$ is isomorphic to the cube graph.

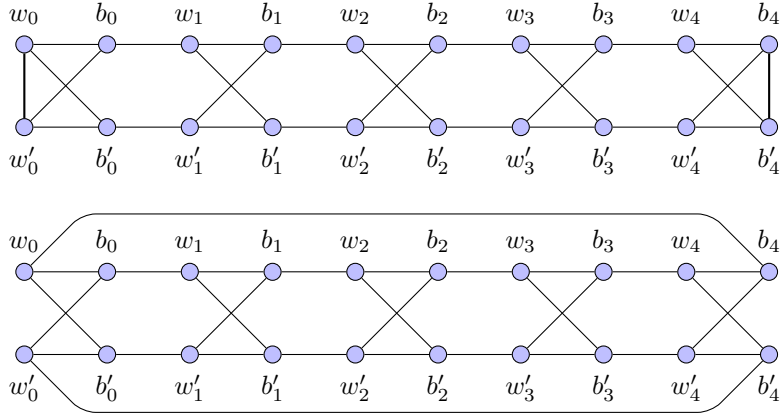
\begin{figure}
\begin{center}
\begin{tikzpicture}[scale=1.1]
\tikzstyle{vertex}=[fill=blue!25, circle, draw=black, inner sep=0.8mm]


\foreach \x in {0,1,2,3,4,5,6,7,8,9} {
\ifthenelse{\isodd{\x}}{ 
\pgfmathsetmacro\result{(\x -1)/2}
\node [vertex, label=below:{\footnotesize $b_{\pgfmathprintnumber{\result}}'$}] (v\x) at (\x,0) {};
\node [vertex, label=above:{\footnotesize $b_{\pgfmathprintnumber{\result}}$}] (w\x) at (\x,1) {}; }
{
\pgfmathsetmacro\result{\x /2 }
\node [vertex, label=below:{\footnotesize $w_{\pgfmathprintnumber{\result}}'$}] (v\x) at (\x,0) {};
\node [vertex, label=above:{\footnotesize $w_{\pgfmathprintnumber{\result}}$}] (w\x) at (\x,1) {};
}
}

\draw (v0)--(v1)--(v2)--(v3)--(v4)--(v5)--(v6)--(v7)--(v8)--(v9);
\draw (w0)--(w1)--(w2)--(w3)--(w4)--(w5)--(w6)--(w7)--(w8)--(w9);
\draw [thick] (v0)--(w0);
\draw (v0)--(w1);
\draw (v1)--(w0);
\draw (v2)--(w3);
\draw (v3)--(w2);
\draw (v4)--(w5);
\draw (v5)--(w4);
\draw (v6)--(w7);
\draw (v7)--(w6);
\draw (v8)--(w9);
\draw (v9)--(w8);
\draw [thick] (v9)--(w9);
\pgftransformyshift{-2.75cm}
\foreach \x in {0,1,2,3,4,5,6,7,8,9} {
\ifthenelse{\isodd{\x}}{ 
\pgfmathsetmacro\result{(\x -1)/2}
\node [vertex, label=below:{\footnotesize $b_{\pgfmathprintnumber{\result}}'$}] (v\x) at (\x,0) {};
\node [vertex, label=above:{\footnotesize $b_{\pgfmathprintnumber{\result}}$}] (w\x) at (\x,1) {}; }
{
\pgfmathsetmacro\result{\x /2 }
\node [vertex, label=below:{\footnotesize $w_{\pgfmathprintnumber{\result}}'$}] (v\x) at (\x,0) {};
\node [vertex, label=above:{\footnotesize $w_{\pgfmathprintnumber{\result}}$}] (w\x) at (\x,1) {};
}
}
\draw (v0)--(v1)--(v2)--(v3)--(v4)--(v5)--(v6)--(v7)--(v8)--(v9);
\draw (w0)--(w1)--(w2)--(w3)--(w4)--(w5)--(w6)--(w7)--(w8)--(w9);
\draw (v0)--(w1);
\draw (v1)--(w0);
\draw (v2)--(w3);
\draw (v3)--(w2);
\draw (v4)--(w5);
\draw (v5)--(w4);
\draw (v6)--(w7);
\draw (v7)--(w6);
\draw (v8)--(w9);
\draw (v9)--(w8);

\draw [rounded corners = 5pt] (v0) -- (0.7,-0.7) -- (8.3,-0.7) -- (v9);
\draw [rounded corners = 5pt] (w0) -- (0.7,1.7) -- (8.3,1.7)--(w9);


\end{tikzpicture}
\caption{The \ks{} graph (top) and \gm{} graph (bottom)}
\label{ksfig}
\end{center}
\end{figure}


The \emph{bipartite double} of a graph $G$ is the tensor product $G \times K_2$ with vertex set $V(G) \times \{0,1\}$ and where 
\[
\{(u,i), (v,j)\} \in E(G \times K_2) \Longleftrightarrow (u,v) \in E(G) \text{  and  } i \ne j. 
\]
In other words, each vertex of $G$ is replaced by a pair of non-adjacent vertices and each edge of $G$ is replaced by a ``crossed matching'' between the corresponding vertex-pairs. If $G$ is connected and not bipartite, then $G \times K_2$ is a connected bipartite graph. If $G$ is bipartite, then $G \times K_2$ is just two copies of $G$. 

If $v$ is an eigenvector for $G$ with eigenvalue $\lambda$, then we can form two related eigenvalues for $G \times K_2$ with eigenvalues $\pm \lambda$. One is obtained by assigning the value $v_x$ to each of the vertices $(x,0)$ and $(x,1)$ and the other is obtained by assigning the value $v_x$ to the vertex $(x,0)$ and $-v_x$ to the vertex $(x,1)$.

\begin{lemma} \label{bipart}
If $G$ is a graph with spectrum $\text{sp}(G)$, then 
\[
\mathrm{sp}(G \times K_2) = \mathrm{sp}(G) \uplus -\mathrm{sp(G)},
\]
where $\uplus$ denotes multiset union.
In particular, $G$ has no eigenvalues in $(-1,1)$ if and only if $G \times K_2$ has no eigenvalues in $(-1,1)$. \qed
\end{lemma}

\begin{lemma}
The bipartite double of the \ks{} graph $\KS{k}$  is the \gm{} graph $\GM{2k}$. 
\end{lemma}

\begin{proof}
Because $B(k)$ is bipartite, $B(k) \times K_2$ consists of two disjoint copies of $B(k)$. To obtain the bipartite double of the {\ks} graph, we need to add the matchings coming from the original edges $w_0w_0'$ and $b_{k-1}b_{k-1}'$. This connects both ends of the two copies of $B(k)$ resulting in the graph $\GM{2k}$.
\end{proof}

Thus half of the \gm{} graphs---those where the number of vertices is a multiple of $8$---are bipartite doubles of the \ks{} graphs. We will see in \cref{thm:KS-GM-bip-double} that, while the other \gm{} graphs are bipartite, they are not the bipartite double of a non-bipartite graph. 

\subsection{Sporadic graphs}

\cref{tab:maintable} lists all the cubic graphs with no eigenvalues in $(-1,1)$ that are not in the two infinite families, also identifying those that are bipartite. They are shown in Figure \ref{fig:allexamples}, in the same order they appear in the table. 

Most of these graphs either have common names already, or are easily described in terms of simple operations on named graphs. In particular, recall that \emph{truncating} a vertex in a cubic graph refers to replacing that vertex with a triangle (also known as a $Y$-$\Delta$ operation).

\begin{table}[htbp]
\begin{center}
\begin{tabular}{ccll}
\toprule
Vertices & Bipartite & \mc{Description} \\
\midrule
$4$ & & Complete graph $K_4$.\\
$10$ & & Petersen graph. \\
$10$ & & $K_{3,3}$ with two white vertices truncated. \\
$12$ & & Petersen with one vertex truncated. \\
$12$ & & $K_{3,3}$ with three white vertices truncated.\\
$14$ &\checkmark & Heawood graph.\\
$16$ & \checkmark & M\"obius-Kantor graph a.k.a generalized Petersen $G(8,3)$.\\
$16$ & & Cube with four white vertices truncated.\\
$20$ & \checkmark & Desargues graph a.k.a. generalized Petersen $G(10,3)$.\\
$20$ & \checkmark & Cospectral mate for Desargues graph.\\
$24$ & \checkmark & Bipartite double (of either $12$-vertex example).\\
$24$ & \checkmark & Bicirculant (see row $4$, column $3$ of Figure~\ref{fig:allexamples}).\\
$32$ & \checkmark & See row $5$, column $1$ of Figure~\ref{fig:allexamples}.\\
$32$ & \checkmark & Bipartite double of cube with four vertices truncated. \\
\bottomrule
\end{tabular}
\end{center}
\caption{Sporadic cubic graphs with no eigenvalues in $(-1,1)$}
\label{tab:maintable}
\end{table}


\input{allexamplesfigures}

As a side note, an astute reader might note that the cubic  graphs with integer eigenvalues, which have been classified in \cite{BusCve1976}, that do not have  $0$ as an eigenvalue should appear in this list and indeed they do. Of the $13$ cubic graphs with integer eigenvalues, there are $6$ without $0$ as an eigenvalue: $K_4$, the cube, the Petersen graph, $K_{3,3}$ with two white vertices truncated, the Desargues graph and its cospectral mate.

%% file: allexamplesfigures.tex
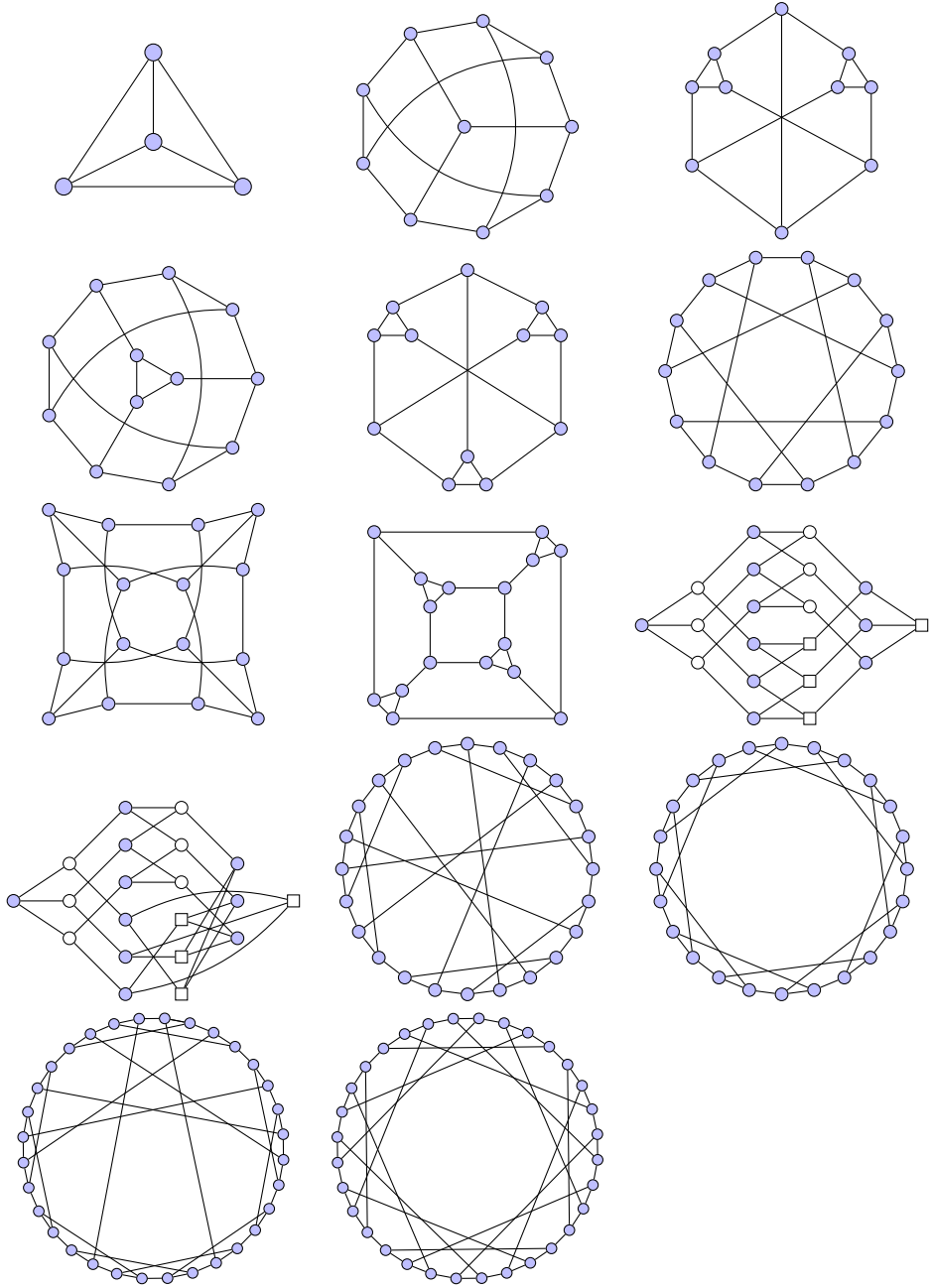
\begin{figure}[p]
\centering
\begin{tabular}{ccc}
\begin{tikzpicture}[scale=0.6]
\tikzstyle{vertex}=[fill=blue!25, circle, draw=black, inner sep=0.8mm]
		\node [vertex] (0) at (0, 2) {};
		\node [vertex] (1) at (-2, -1) {};
		\node [vertex] (2) at (0, 0) {};
		\node [vertex] (3) at (2, -1) {};
 \node (4) at (0,-1.8) {\phantom{aa}};
		\draw (0) to (2);
		\draw (2) to (1);
		\draw (1) to (3);
		\draw (3) to (2);
		\draw (0) to (1);
		\draw (0) to (3);
\end{tikzpicture}     & 
\begin{tikzpicture}[scale=0.6]
\tikzstyle{vertex}=[fill=blue!25, circle, draw=black, inner sep=0.6mm]
\foreach \x in {0,1,...,8} {
  \node[vertex] (v\x) at (\x*40:2.4cm) {};
}
\node[vertex] (v9) at (0,0) {};
\draw (v0) -- (v1) -- (v2) -- (v3) -- (v4) -- (v5) -- (v6) -- (v7) -- (v8) -- (v0);
\draw (v9) -- (v0);
\draw (v9) -- (v3);
\draw (v9) -- (v6);
\draw (v1) to [bend right] (v5);
\draw (v2) to [bend left] (v7);
\draw (v4) to [bend right]  (v8);
\end{tikzpicture}&
\begin{tikzpicture}[scale=0.6]
\tikzstyle{vertex}=[fill=blue!25, circle, draw=black, inner sep=0.6mm]
		\node [vertex] (0) at (0, 3) {};
		\node [vertex] (1) at (0, -2) {};
		\node [vertex] (2) at (-2, -0.5) {};
		\node [vertex] (3) at (2, -0.5) {};
		\node [vertex] (4) at (-1.25, 1.25) {};
		\node [vertex] (5) at (-2, 1.25) {};
		\node [vertex] (6) at (-1.5, 2) {};
		\node [vertex] (7) at (1.5, 2) {};
		\node [vertex] (8) at (1.25, 1.25) {};
		\node [vertex] (9) at (2, 1.25) {};
		\draw (0) to (6);
		\draw (6) to (5);
		\draw (5) to (2);
		\draw (2) to (1);
		\draw (1) to (3);
		\draw (3) to (9);
		\draw (9) to (7);
		\draw (7) to (0);
		\draw (0) to (1);
		\draw (2) to (8);
		\draw (8) to (7);
		\draw (8) to (9);
		\draw (3) to (4);
		\draw (4) to (5);
		\draw (4) to (6);
  \end{tikzpicture} \\
  \begin{tikzpicture}[scale=0.6]
\tikzstyle{vertex}=[fill=blue!25, circle, draw=black, inner sep=0.6mm]
\foreach \x in {0,1,...,8} {
  \node[vertex] (v\x) at (\x*40:2.4cm) {};
}
\foreach \x in {9,10,11} {
  \node[vertex] (v\x) at (\x*120:0.6cm) {};
}
\draw (v0) -- (v1) -- (v2) -- (v3) -- (v4) -- (v5) -- (v6) -- (v7) -- (v8) -- (v0);
\draw (v9) -- (v10) -- (v11) -- (v9);
\draw (v9) -- (v0);
\draw (v10) -- (v3);
\draw (v11) -- (v6);
\draw (v1) to [bend right] (v5);
\draw (v2) to [bend left] (v7);
\draw (v4) to [bend right]  (v8);
\end{tikzpicture} &
   \begin{tikzpicture}[scale = 0.5]
  \tikzstyle{vertex}=[fill=blue!25, circle, draw=black, inner sep=0.6mm]
		\node [vertex] (1) at (2, 1.75) {};
		\node [vertex] (2) at (0, -2.25) {};
		\node [vertex] (3) at (-2.5, -1.5) {};
		\node [vertex] (4) at (2.5, -1.5) {};
		\node [vertex] (5) at (0, 2.75) {};
		\node [vertex] (6) at (-2.5, 1) {};
		\node [vertex] (7) at (-1.5, 1) {};
		\node [vertex] (9) at (-2, 1.75) {};
		\node [vertex] (10) at (2.5, 1) {};
		\node [vertex] (11) at (1.5, 1) {};
		\node [vertex] (12) at (-0.5, -3) {};
		\node [vertex] (13) at (0.5, -3) {};
		\draw (5) to (2);
		\draw (2) to (12);
		\draw (12) to (13);
		\draw (13) to (2);
		\draw (5) to (9);
		\draw (9) to (6);
		\draw (6) to (7);
		\draw (7) to (9);
		\draw (6) to (3);
		\draw (3) to (12);
		\draw (13) to (4);
		\draw (4) to (10);
		\draw (10) to (11);
		\draw (11) to (1);
		\draw (1) to (10);
		\draw (1) to (5);
		\draw (11) to (3);
		\draw (7) to (4);
\end{tikzpicture}    & 
\begin{tikzpicture}[scale=0.6]
\tikzstyle{vertex}=[fill=blue!25, circle, draw=black, inner sep=0.6mm]
\foreach \x in {0,1,...,13} {
  \node[vertex] (v\x) at (\x*25.71:2.6cm) {};
}
\draw (v0) -- (v1) -- (v2) -- (v3) -- (v4) -- (v5) -- (v6) -- (v7) -- (v8) -- (v9) -- (v10) -- (v11) -- (v12) -- (v13) -- (v0);
\draw (v0) -- (v5);
\draw (v2) -- (v7);
\draw (v4) -- (v9);
\draw (v6) -- (v11);
\draw (v8) -- (v13);
\draw (v10) -- (v1);
\draw (v12) -- (v3);
\end{tikzpicture} \\
\begin{tikzpicture}[scale=0.4]
\tikzstyle{vertex}=[fill=blue!25, circle, draw=black, inner sep=0.6mm]
		\node [vertex] (0) at (-1, 1) {};
		\node [vertex] (1) at (-1, -1) {};
		\node [vertex] (2) at (1, -1) {};
		\node [vertex] (3) at (1, 1) {};
		\node [vertex] (4) at (-1.5, 3) {};
		\node [vertex] (5) at (1.5, 3) {};
		\node [vertex] (6) at (-1.5, -3) {};
		\node [vertex] (7) at (1.5, -3) {};
		\node [vertex] (8) at (-3, 1.5) {};
		\node [vertex] (9) at (-3, -1.5) {};
		\node [vertex] (10) at (3, 1.5) {};
		\node [vertex] (11) at (3, -1.5) {};
		\node [vertex] (12) at (3.5, 3.5) {};
		\node [vertex] (13) at (-3.5, 3.5) {};
		\node [vertex] (14) at (-3.5, -3.5) {};
		\node [vertex] (15) at (3.5, -3.5) {};
		\draw (13) to (8);
		\draw (8) to (9);
		\draw (9) to (14);
		\draw (14) to (6);
		\draw (6) to (7);
		\draw (7) to (15);
		\draw (15) to (11);
		\draw (11) to (10);
		\draw (10) to (12);
		\draw (12) to (5);
		\draw (5) to (4);
		\draw (4) to (13);
		\draw (13) to (0);
		\draw (3) to (12);
		\draw (2) to (15);
		\draw (1) to (14);
		\draw [bend right=15] (9) to (2);
		\draw [bend right=15] (1) to (11);
		\draw [bend left=15] (0) to (10);
		\draw [bend left=15] (8) to (3);
		\draw [bend right=15] (4) to (1);
		\draw [bend right=15] (0) to (6);
		\draw [bend left=15] (3) to (7);
		\draw [bend right=15] (2) to (5);
\end{tikzpicture}&
\begin{tikzpicture}[scale=0.5]
\tikzstyle{vertex}=[fill=blue!25, circle, draw=black, inner sep=0.6mm]
		\node [vertex] (1) at (1, 1) {};
		\node [vertex] (2) at (-1, -1) {};
		\node [vertex] (4) at (-2.5, 2.5) {};
		\node [vertex] (5) at (1.75, 1.75) {};
		\node [vertex] (6) at (2.5, -2.5) {};
		\node [vertex] (8) at (-1.25, 1.25) {};
		\node [vertex] (9) at (-0.5, 1) {};
		\node [vertex] (10) at (-1, 0.5) {};
		\node [vertex] (11) at (0.5, -1) {};
		\node [vertex] (12) at (1, -0.5) {};
		\node [vertex] (13) at (1.25, -1.25) {};
		\node [vertex] (14) at (-1.75, -1.75) {};
		\node [vertex] (15) at (-2.5, -2) {};
		\node [vertex] (16) at (-2, -2.5) {};
		\node [vertex] (17) at (2, 2.5) {};
		\node [vertex] (18) at (2.5, 2) {};
		\draw (4) to (17);
		\draw (17) to (18);
		\draw (18) to (6);
		\draw (6) to (16);
		\draw (16) to (15);
		\draw (15) to (4);
		\draw (4) to (8);
		\draw (8) to (10);
		\draw (10) to (9);
		\draw (9) to (8);
		\draw (10) to (2);
		\draw (2) to (14);
		\draw (14) to (15);
		\draw (14) to (16);
		\draw (2) to (11);
		\draw (11) to (12);
		\draw (12) to (13);
		\draw (13) to (11);
		\draw (13) to (6);
		\draw (12) to (1);
		\draw (1) to (9);
		\draw (1) to (5);
		\draw (5) to (17);
		\draw (5) to (18);
\end{tikzpicture} &
 \begin{tikzpicture}[scale=0.5]
\tikzstyle{vertex}=[fill=blue!25, circle, draw=black, inner sep=0.6mm]
\tikzstyle{vertex2}=[fill=white, circle, draw=black, inner sep=0.6mm]
\tikzstyle{vertex3}=[fill=white, rectangle, draw=black, inner sep=0.8mm]
		\node [vertex] (0) at (-1, 0.5) {};
		\node [vertex] (1) at (-1, -0.5) {};
		\node [vertex] (2) at (-1, -1.5) {};
		\node [vertex] (3) at (-1, -2.5) {};
		\node [vertex3] (4) at (0.5, -2.5) {};
		\node [vertex3] (5) at (0.5, -1.5) {};
		\node [vertex3] (6) at (0.5, -0.5) {};
		\node [vertex2] (7) at (0.5, 0.5) {};
		\node [vertex2] (8) at (0.5, 1.5) {};
		\node [vertex2] (9) at (0.5, 2.5) {};
		\node [vertex] (10) at (-1, 2.5) {};
		\node [vertex] (11) at (-1, 1.5) {};
		\node [vertex2] (12) at (-2.5, 0) {};
		\node [vertex2] (13) at (-2.5, -1) {};
		\node [vertex2] (14) at (-2.5, 1) {};
		\node [vertex] (15) at (2, 0) {};
		\node [vertex] (16) at (2, 1) {};
		\node [vertex] (17) at (2, -1) {};
		\node [vertex] (18) at (-4, 0) {};
		\node [vertex3] (19) at (3.5, 0) {};
		\draw (18) to (14);
		\draw (14) to (10);
		\draw (18) to (12);
		\draw (18) to (13);
		\draw (13) to (3);
		\draw (19) to (16);
		\draw (19) to (15);
		\draw (19) to (17);
		\draw (17) to (4);
		\draw (16) to (9);
		\draw (10) to (9);
		\draw (9) to (11);
		\draw (11) to (7);
		\draw (7) to (0);
		\draw (0) to (8);
		\draw (8) to (10);
		\draw (1) to (6);
		\draw (6) to (2);
		\draw (2) to (4);
		\draw (4) to (3);
		\draw (3) to (5);
		\draw (5) to (1);
		\draw (12) to (11);
		\draw (14) to (1);
		\draw (13) to (0);
		\draw (12) to (2);
		\draw (5) to (15);
		\draw (15) to (8);
		\draw (6) to (16);
		\draw (7) to (17);
\end{tikzpicture}    \\
\begin{tikzpicture}[scale=0.5]
\tikzstyle{vertex}=[fill=blue!25, circle, draw=black, inner sep=0.6mm]
\tikzstyle{vertex2}=[fill=white, circle, draw=black, inner sep=0.6mm]
\tikzstyle{vertex3}=[fill=white, rectangle, draw=black, inner sep=0.8mm]
		\node [vertex] (0) at (-1, 0.5) {};
		\node [vertex] (1) at (-1, -0.5) {};
		\node [vertex] (2) at (-1, -1.5) {};
		\node [vertex] (3) at (-1, -2.5) {};
		\node [vertex3] (4) at (0.5, -2.5) {};
		\node [vertex3] (5) at (0.5, -1.5) {};
		\node [vertex3] (6) at (0.5, -0.5) {};
		\node [vertex2] (7) at (0.5, 0.5) {};
		\node [vertex2] (8) at (0.5, 1.5) {};
		\node [vertex2] (9) at (0.5, 2.5) {};
		\node [vertex] (10) at (-1, 2.5) {};
		\node [vertex] (11) at (-1, 1.5) {};
		\node [vertex2] (12) at (-2.5, 0) {};
		\node [vertex2] (13) at (-2.5, -1) {};
		\node [vertex2] (14) at (-2.5, 1) {};
		\node [vertex] (15) at (2, 0) {};
		\node [vertex] (16) at (2, 1) {};
		\node [vertex] (17) at (2, -1) {};
		\node [vertex] (18) at (-4, 0) {};
		\node [vertex3] (19) at (3.5, 0) {};
		\draw (18) to (14);
		\draw (14) to (10);
		\draw (18) to (12);
		\draw (18) to (13);
		\draw (13) to (3);
		\draw (5) to (16);
		\draw (6) to (15);
		\draw (6) to (17);
		\draw (17) to (5);
		\draw (16) to (9);
		\draw (10) to (9);
		\draw (9) to (11);
		\draw (11) to (7);
		\draw (7) to (0);
		\draw (0) to (8);
		\draw (8) to (10);
		\draw (1) to [bend left = 20] (19);
		\draw (19) to (2);
		\draw (2) to (5);
		\draw (6) to (3);
		\draw (3) to [bend right = 18] (19);
		\draw (4) to (1);
		\draw (12) to (11);
		\draw (14) to (1);
		\draw (13) to (0);
		\draw (12) to (2);
		\draw (4) to (15);
		\draw (15) to (8);
		\draw (4) to (16);
		\draw (7) to (17);
\end{tikzpicture} & 
\begin{tikzpicture}[scale=0.7]
\tikzstyle{vertex}=[fill=blue!25, circle, draw=black, inner sep=0.6mm]
\foreach \x in {0,1,...,23} {
  \node[vertex] (v\x) at (\x*15:2.4cm) {};
}
\draw (v0)--(v1)--(v2)--(v3)--(v4)--(v5)--(v6)--(v7)--(v8)--(v9)--(v10)--(v11)--(v12)--(v13)--(v14)--(v15)--(v16)--(v17)--(v18)--(v19)--(v20)--(v21)--(v22)--(v23)--(v0);
\draw (v0)--(v5);
\draw (v1)--(v12);
\draw (v2)--(v7);
\draw (v3)--(v14);
\draw (v4)--(v17);
\draw (v6)--(v19);
\draw (v8)--(v13);
\draw (v9)--(v20);
\draw (v10)--(v15);
\draw (v11)--(v22);
\draw (v16)--(v21);
\draw (v18)--(v23);
\end{tikzpicture} &
\begin{tikzpicture}[scale=0.7]
\tikzstyle{vertex}=[fill=blue!25, circle, draw=black, inner sep=0.6mm]
\foreach \x in {0,1,...,23} {
  \node[vertex] (v\x) at (\x*15:2.4cm) {};
}
\draw (v0)--(v1)--(v2)--(v3)--(v4)--(v5)--(v6)--(v7)--(v8)--(v9)--(v10)--(v11)--(v12)--(v13)--(v14)--(v15)--(v16)--(v17)--(v18)--(v19)--(v20)--(v21)--(v22)--(v23)--(v0);
\draw (v0)--(v5);
\draw (v2)--(v7);
\draw (v4)--(v9);
\draw (v6)--(v11);
\draw (v8)--(v13);
\draw (v10)--(v15);
\draw (v12)--(v17);
\draw (v14)--(v19);
\draw (v16)--(v21);
\draw (v18)--(v23);
\draw (v20)--(v1);
\draw (v22)--(v3);
\end{tikzpicture} \\
\begin{tikzpicture}[scale=0.7]
\tikzstyle{vertex}=[fill=blue!25, circle, draw=black, inner sep=0.5mm]
\foreach \x in {0,1,...,31} {
  \node[vertex] (v\x) at (85+\x*360/32:2.5cm) {};
}
\draw (v0)--(v1)--(v2)--(v3)--(v4)--(v5)--(v6)--(v7)--(v8)--(v9)--(v10)--(v11)--(v12)--(v13)--(v14)--(v15)--(v16)--(v17)--(v18)--(v19)--(v20)--(v21)--(v22)--(v23)--(v24)--(v25)--(v26)--(v27)--(v28)--(v29)--(v30)--(v31)--(v0);
\draw (v2)--(v29);
\draw (v0)--(v31);
\draw (v0)--(v19);
\draw (v5)--(v10);
\draw (v3)--(v24);
\draw (v9)--(v30);
\draw (v4)--(v31);
\draw (v13)--(v18);
\draw (v6)--(v25);
\draw (v7)--(v12);
\draw (v11)--(v16);
\draw (v1)--(v14);
\draw (v8)--(v27);
\draw (v21)--(v26);
\draw (v23)--(v28);
\draw (v17)--(v22);
\draw (v15)--(v20);
\end{tikzpicture}& 
\begin{tikzpicture}[scale=0.7]
\tikzstyle{vertex}=[fill=blue!25, circle, draw=black, inner sep=0.5mm]
\foreach \x in {0,1,...,31} {
  \node[vertex] (v\x) at (85+\x*360/32:2.5cm) {};
}
\draw (v0)--(v1)--(v2)--(v3)--(v4)--(v5)--(v6)--(v7)--(v8)--(v9)--(v10)--(v11)--(v12)--(v13)--(v14)--(v15)--(v16)--(v17)--(v18)--(v19)--(v20)--(v21)--(v22)--(v23)--(v24)--(v25)--(v26)--(v27)--(v28)--(v29)--(v30)--(v31)--(v0);
\draw (v0)--(v9);
\draw (v1)--(v24);
\draw (v2)--(v11);
\draw (v3)--(v26);
\draw (v4)--(v29);
\draw (v5)--(v12);
\draw (v6)--(v15);
\draw (v7)--(v30);
\draw (v8)--(v17);
\draw (v10)--(v19);
\draw (v13)--(v20);
\draw (v14)--(v23);
\draw (v16)--(v25);
\draw (v18)--(v27);
\draw (v21)--(v28);
\draw (v22)--(v31);
\end{tikzpicture}& 
\end{tabular}
\caption{Sporadic graphs with no eigenvalues in $(-1,1)$. \textit{From left to right.} \textit{Row 1: } complete graph on $4$ vertices; Petersen graph;  $K_{3,3}$ with two non-adjacent vertices truncated. \textit{Row 2: } Petersen graph with one vertex truncated;   $K_{3,3}$ with three pairwise non-adjacent vertices 
    truncated; Heawood graph. \textit{Row 3: } M\"{o}bius-Kantor graph; cube with one bipartite colour class truncated; Desargues graph. \textit{Row 4: } cospectral mate of Desargues graph (obtained from Desargues graph by Godsil-McKay switching on the square vertices); bipartite double of Petersen graph with one vertex truncated, bicirculant graph on $24$ vertices.  \textit{Row 5:} graph on $32$ vertices; bipartite double of the cube with one bipartite colour class truncated. \label{fig:allexamples}}
\end{figure}

%% file: proof.tex
\section{Proofs}

We start with an elementary but very useful observation:

\begin{lemma}
A graph $G$ has no eigenvalues in $(-1,1)$ if and only if the matrix $M(G) = A(G)^2 - I$ is positive semi-definite. \qed
\end{lemma}

The reason this is so useful is because a matrix is positive semi-definite if and only if all of its principal minors are non-negative. If a particular local configuration in $G$ would create a negative principal minor of $M(G)$, then that configuration is forbidden in any graph with no eigenvalues in $(-1,1)$. 

The proof strategy will be to first classify the \emph{bipartite} cubic graphs with no eigenvalues in $(-1,1)$, and then to determine which of those graphs are bipartite doubles of non-bipartite graphs.

\begin{theorem}
Let $G$ be a cubic bipartite graph with no eigenvalues in $(-1,1)$. Then the girth of $G$ is either $4$ or $6$. 
Furthermore, if $G$ has girth $4$, then it is a Guo-Mohar graph, and if $G$ has girth $6$, then it is one of the eight graphs indicated as bipartite in Table~\ref{tab:maintable}.
\end{theorem}

We divide the proof of this into two cases, according to whether the girth is equal to $4$ or greater than $4$, with the next two subsections dealing with these cases in turn.


\subsection{Bipartite with girth equal to \texorpdfstring{$4$}{4}}\label{bip4}

\input girth4

\subsection{Bipartite with girth at least \texorpdfstring{$6$}{6}}\label{bip6}

\input girth6

\subsection{Non-bipartite graphs}\label{nonbip}

\input nonbipartite

\input maximalgapset

%% file: girth4.tex
In this section, we assume that $G$ is a cubic bipartite graph of girth $4$ with no eigenvalues in $(-1,1)$ and show that $G$ is necessarily a Guo-Mohar graph. 

\begin{figure}[htb]
\begin{center}
\begin{tikzpicture}[scale=1.25]
\tikzstyle{bvertex}=[fill=black, circle, draw=black, inner sep=0.8mm]
\tikzstyle{wvertex}=[fill=white, circle, draw=black, inner sep=0.8mm]
\foreach \x in {0,2,4,6,8,10} {
  \node [bvertex] (v\x) at (15+30*\x:1.75cm) {};
  \node [bvertex] (w\x) at (15+30*\x:2.5cm) {};
}
\foreach \x in {1,3,5,7,9,11} {
  \node [wvertex] (v\x) at (15+30*\x:1.75cm) {};
  \node [wvertex] (w\x) at (15+30*\x:2.5cm) {};
}

\draw (v0)--(v1)--(v2)--(v3)--(v4)--(v5)--(v6)--(v7)--(v8)--(v9)--(v10)--(v11)--(v0);
\draw (w0)--(w1)--(w2)--(w3)--(w4)--(w5)--(w6)--(w7)--(w8)--(w9)--(w10)--(w11)--(w0);
\draw (v0)--(w1);
\draw (w0)--(v1);
\draw (v2)--(w3);
\draw (w2)--(v3);
\draw (v4)--(w5);
\draw (w4)--(v5);
\draw (v6)--(w7);
\draw (w6)--(v7);
\draw (v8)--(w9);
\draw (w8)--(v9);
\draw (v10)--(w11);
\draw (w10)--(v11);
\foreach \x in {0,1,...,5} {
  \node at (90-60*\x:2.75cm) {$C_\x$};
}
\end{tikzpicture}
\end{center}
\caption{Guo-Mohar graph $\mathrm{GM}(6)$}
\label{fig:GuoMohar6}
\end{figure}
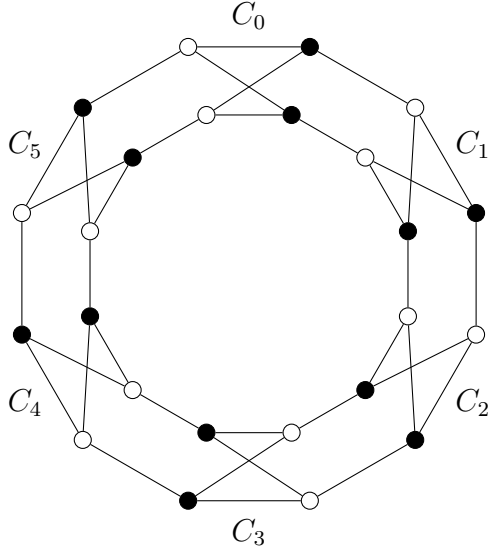

Recall that Guo-Mohar graph $\GM{k}$ consists of $k>1$ induced $4$-cycles $C_0$, $C_1$, \ldots, $C_{k-1}$ connected in a cyclic order; see  \cref{fig:GuoMohar6} for an example. We use the notation established in \cref{families}, where each $4$-cycle $C_i = (w_i, b_i, w'_i, b'_i)$ and the cycle $C_i$ is connected to $C_{i+1}$ by the two edges $b_iw_{i+1}$ and $b'_iw'_{i+1}$ (with all indices being taken modulo $k$), as illustrated in  \cref{ksfig}.

We shall show that if $G$ contains a subgraph isomorphic to the base graph $B(k)$ (that is, a \emph{path} of $k$ induced $4$-cycles connected as described above) then either $G$ is isomorphic to $\mathrm{GM}(k)$ or it contains a subgraph $B(k+1)$. We start with the case $k=1$, where the case analysis is a little more complicated than for the case $k>1$.

\begin{lemma}
Let $G$ be a cubic bipartite graph with no eigenvalues in $(-1,1)$ and suppose that $G$ contains an induced $4$-cycle $C_0$. Then $G$ contains a subgraph isomorphic to $B(2)$.
\end{lemma}

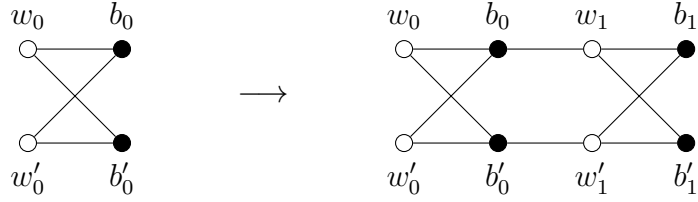
\begin{figure}
\begin{center}
\begin{tikzpicture}[scale=1.25]
\tikzstyle{bvertex}=[fill=black, circle, draw=black, inner sep=0.8mm]
\tikzstyle{wvertex}=[fill=white, circle, draw=black, inner sep=0.8mm]

\node [wvertex, label = below:$w'_0$] (w0) at (0,0) {};
\node [bvertex, label = below:$b'_0$] (b0) at (1,0) {};
\node [wvertex, label = above:$w_0$] (wp0) at (0,1) {};
\node [bvertex, label = above:$b_0$] (bp0) at (1,1) {};

\draw (w0)--(b0)--(wp0)--(bp0)--(w0);

\node at (2.5,0.5) {$\longrightarrow$};

\pgftransformxshift{4cm}

\node [wvertex, label = below:$w'_0$] (w0) at (0,0) {};
\node [bvertex, label = below:$b'_0$] (b0) at (1,0) {};
\node [wvertex, label = above:$w_0$] (wp0) at (0,1) {};
\node [bvertex, label = above:$b_0$] (bp0) at (1,1) {};

\draw (w0)--(b0)--(wp0)--(bp0)--(w0);

\node [wvertex, label = below:$w'_1$] (w1) at (2,0) {};
\node [bvertex, label = below:$b'_1$] (b1) at (3,0) {};
\node [wvertex, label = above:$w_1$] (wp1) at (2,1) {};
\node [bvertex, label = above:$b_1$] (bp1) at (3,1) {};

\draw (w1)--(b1)--(wp1)--(bp1)--(w1);

\draw (bp0)--(wp1);
\draw (b0)--(w1);

\end{tikzpicture}
\end{center}
\caption{One $4$-cycle leads to another}
\label{fig:casek=1}
\end{figure}

\begin{proof}
Starting with $C_0$, we know that $b_0$ has a third neighbour, say $w_1$. Now $w_1$ is not adjacent to $b'_0$ because otherwise $b_0$ and $b'_0$ would be twin vertices and $A(G)$ would have repeated rows and hence an eigenvalue of zero. So the third neighbour of $b'_0$, say $w'_1$, is distinct from $w_1$.

Now $w_1$ in turn is adjacent to two new vertices, $b_1$, $b'_1$; since the graph is bipartite, $w_1$ is not adjacent to either of $w_0,w_0'$ and $w_1$ is not adjacent to $b_0'$, which already has three distinct neighbours. The vertex $b_1$ is not adjacent to both of the vertices $w_0$, $w'_0$, because, if it were, then $\{w_0 ,w'_0\}$ would be a pair of twin vertices. See \cref{fig:casek=1a} for our current state.

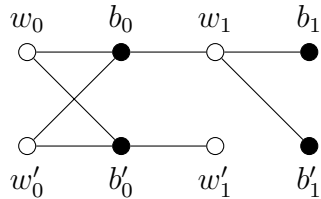
\begin{figure}[htb]
\begin{center}
\begin{tikzpicture}[scale=1.25]
\tikzstyle{bvertex}=[fill=black, circle, draw=black, inner sep=0.8mm]
\tikzstyle{wvertex}=[fill=white, circle, draw=black, inner sep=0.8mm]

\pgftransformxshift{4cm}

\node [wvertex, label = below:$w'_0$] (w0) at (0,0) {};
\node [bvertex, label = below:$b'_0$] (b0) at (1,0) {};
\node [wvertex, label = above:$w_0$] (wp0) at (0,1) {};
\node [bvertex, label = above:$b_0$] (bp0) at (1,1) {};

\draw (w0)--(b0)--(wp0)--(bp0)--(w0);

\node [wvertex, label = below:$w'_1$] (w1) at (2,0) {};
\node [bvertex, label = below:$b'_1$] (b1) at (3,0) {};
\node [wvertex, label = above:$w_1$] (wp1) at (2,1) {};
\node [bvertex, label = above:$b_1$] (bp1) at (3,1) {};

\draw (wp1)--(b1);
\draw (bp1)--(wp1);

\draw (bp0)--(wp1);
\draw (b0)--(w1);
\end{tikzpicture}
\end{center}
\caption{The $4$-cycle, $\{w_0,b_0,w_0',b_0'\}$, shown with other neighbours, $w_1,w_1'$, of $b_0$ and $b_0'$, respectively, and the neighbours of $w_1$.}
\label{fig:casek=1a}
\end{figure}

We aim to show that $b_1$ is adjacent to $w'_1$, so assume that this is not the case. Then if we take $S = \{b_0, b'_0, b_1\}$ we have
\[
M_{SS} = 
\left[
\begin{array}{ccc}
2 & 2 & 1+x\\
2 & 2 & x \\
1+x & x & 2
\end{array}
\right]
\]
where $x$ is either $1$ or $0$, when  $b_1$ is adjacent to a vertex in $\{w_0, w'_0\}$ or not, respectively. 

It is easy to verify that $\det M_{SS} = -2$,  independent of the value of $x$, and thus $M$ is not positive semi-definite. Therefore our assumption was incorrect and it follows that $b_1$ is adjacent to $w'_1$.  A symmetric argument shows that $b'_1$ is adjacent to $w_1$. Therefore $G$ contains a configuration of two induced $4$-cycles forming the graph $B(2)$, as depicted in the second diagram of \cref{fig:casek=1}.
\end{proof}

The next lemma shows that if there is any edge from the last cycle in the path to the first, then $G$ is a Guo-Mohar graph.

\begin{lemma}
Let $G$ be a cubic bipartite graph with no eigenvalues in $(-1,1)$. Suppose that $G$ contains a subgraph $B(k)$ for some $k>1$. If either of the black vertices $\{b_{k-1}, b'_{k-1}\}$ is adjacent to either of the white vertices $\{w_0, w'_0\}$, then $G$ is the Guo-Mohar graph $\mathrm{GM}(k)$.
\end{lemma}

\begin{proof}
Without loss of generality, suppose that $b_{k-1}$ is adjacent to $w_0$. Then it is not adjacent to $w'_0$, for that would create twin vertices $\{w_0, w'_0\}$.  If $b'_{k-1}$ is not adjacent to $w'_0$ then let $S = \{w_0, w'_0, w_{k-1}\}$, and count $2$-paths between these vertices to determine that
\[
M_{SS} = 
\left[
\begin{array}{ccc}
2 & 2 & 2\\
2 & 2 & 1 \\
2 & 1 & 2
\end{array}
\right]
\]
if $k=2$ and 
\[
M_{SS} = 
\left[
\begin{array}{ccc}
2 & 2 & 1\\
2 & 2 & 0 \\
1 & 0 & 2
\end{array}
\right]
\]
if $k > 2$. (These matrices are different because $w_0 \sim b_0 \sim w_1$ is a path from $w_0$ to $w_{k-1}$ only when $k=2$.) Both of these matrices have determinant $-2$ and so $M(G)$ is not positive semidefinite.

Therefore $b'_{k-1}$ is adjacent to $w'_0$ and the graph is cubic and isomorphic to $\mathrm{GM}(k)$. 
\end{proof}

\begin{lemma}\label{leminduct}
Let $G$ be a cubic bipartite graph with no eigenvalues in $(-1,1)$. Suppose that $G$ contains a subgraph isomorphic to $B(k)$ but that $G$ is not isomorphic to $\GM{k}$. Then $G$ contains a subgraph isomorphic $B(k+1)$.
\end{lemma}

\begin{proof}
By the assumption that $G$ is not isomorphic to $\GM{k}$ and the previous lemma, $b_{k-1}$ is not adjacent to either of $\{w_0, w'_0\}$ and so its third neighbour is a new vertex $w_k$. Similarly $b'_{k-1}$ is adjacent to a third vertex, call it $w'_k$ where, to avoid twin vertices, it must be the case that $w_k \not= w'_k$. Now $w_k$ is not adjacent to either of the vertices $\{w_0, w'_0\}$ because $G$ is bipartite. Therefore it has two additional neighbours, which we call $b_k$, $b'_k$. If $b_k$ is not adjacent to $w'_k$ then taking  $S = \{b_{k-1},b'_{k-1},b_k\}$ we get that
\[
M_{SS} = 
\left[
\begin{array}{ccc}
2 & 2 & 1\\
2 & 2 & 0 \\
1 & 0 & 2
\end{array}
\right]
\]
which has determinant $-2$. Therefore $b_k$ is adjacent to $w'_k$ and similarly for $b'_k$ thereby adding an additional induced $4$-cycle to the path of $4$-cycles and creating a subgraph isomorphic to $B(k+1)$.
\end{proof}

\begin{theorem}
If $G$ is a cubic bipartite graph of girth $4$ with no eigenvalues in $(-1,1)$, then $G$ is a Guo-Mohar graph $\mathrm{GM}(k)$ for some $k$.
\end{theorem}

\begin{proof}
Let $k$ be the largest integer such that $G$ contains a subgraph isomorphic to $B(k)$; by \cref{fig:casek=1} we know that $k$ exists and that $k \geq 2$. Then by \cref{leminduct}, $G$ is isomorphic to $\GM{k}$.
\end{proof}

%% file: girth6.tex


For a graph with girth at least $5$, there is at most one path of length $2$ between two vertices and so the non-diagonal entries of $A(G)^2$ are either 0 or 1.

Therefore
\[
A(G)^2-I = 2I + A(D_2(G))
\]
where $D_2(G)$ is the ``$2$-distance graph'' of $G$, in other words the graph on the same vertex set as $G$ where two vertices are adjacent if and only if they have distance two in $G$. 

In particular $A(G)^2-I$ is positive semi-definite if and only if $D_2(G)$ has minimum eigenvalue at least $-2$. One of the seminal results of algebraic graph theory is the characterisation of graphs with least eigenvalue $-2$ by Cameron, Goethals, Seidel and Shult in 1976. Our strategy then will be to use this characterisation to determine all of the possibilities for $D_2(G)$ and hence $G$.

Because $G$ is a cubic bipartite graph of girth at least $6$, the graph $D_2(G)$ is a $6$-regular graph with two connected components. 

\begin{theorem}\cite{CamGoeSeiShu1976} \label{min2}
    If $X$ is a regular connected graph with minimum eigenvalue at least $-2$, then one of the following holds:
    \begin{enumerate}[itemsep=0pt]
        \item $X$ is a line graph, or
        \item $X$ is a cocktail party graph, i.e., $X = \overline{mK_2}$, or
        \item $X$ is an exceptional graph with a representation in the $E_8$ root lattice.
    \end{enumerate}
\end{theorem}

The next three lemmas examine each of these three cases in turn, determining all the possibilities for the components of $D_2(G)$ and therefore all the possibilities for $G$.

Since $G$ is bipartite and has girth at least $6$,  we can view it as the point-line incidence graph of a
partial linear space. As $G$ is a cubic graph, each line of this geometry contains three points,
and each point lies in three lines (i.e., the geometry is $3$-regular). The two components of
$D_2(G)$ are the \emph{collinearity} graph on points and the \emph{concurrence} graph on lines. Because
each line determines a triangle of the collinearity graph (and dually each point determines a
triangle in the concurrence graph), the edges of each component of $D_2(G)$ can be partitioned into triangles. A partition of the edges of a graph into edge-disjoint triangles is called a
\emph{triangle decomposition} of the graph.

\begin{lemma}
Let $G$ be a cubic bipartite graph of girth at least $6$ such that $D_2(G)$ has a connected component isomorphic to a cocktail party graph. Then $G$ is the M\"obius-Kantor graph.
\end{lemma}
\begin{proof}
The only $6$-regular cocktail party graph is $\overline{4K_2}$ which has $8$ vertices and so $G$ has $16$ vertices. Of the $38$ cubic bipartite graphs on $16$ vertices, routine computation shows that the only one of girth at least $6$ is the M\"obius-Kantor graph.
This indeed has no eigenvalues in $(-1,1)$ and has a 2-distance graph with components isomorphic to $\overline{4K_2}$. 
\end{proof}

\begin{lemma}
Let $G$ be a cubic bipartite graph of girth at least $6$ such that $D_2(G)$ has a connected component isomorphic to one of the exceptional graphs of \cref{min2}. Then $G$ has $24$ or $32$ vertices and is one of the four graphs of those orders listed in \cref{tab:maintable}.
\end{lemma}
\begin{proof}
The regular graphs (of any degree) with minimum eigenvalue at least $-2$ that are not linegraphs are explicitly listed in \cite[Table 9.1]{BusCveSei1976}. The only $6$-regular graphs in these lists are five graphs on $12$ vertices and $35$ graphs on $16$ vertices. It is straightforward to verify that just two of the $12$-vertex graphs and two of the $16$-vertex graphs have triangle decompositions.  
\end{proof}

\begin{figure}
\begin{center}
\begin{tikzpicture}
\tikzstyle{vertex}=[fill=white, circle, draw=black, thick, inner sep=1mm]
\foreach \x in {0,1,...,11} {
\node [vertex] (v\x) at (30*\x:2.5cm) {};
}
\draw [ultra thick,green!80!black] (v0)--(v1)--(v3)--(v0);
\draw [ultra thick,green!80!black] (v3)--(v4)--(v6)--(v3);
\draw [ultra thick,green!80!black] (v6)--(v7)--(v9)--(v6);
\draw [ultra thick,green!80!black] (v9)--(v10)--(v0)--(v9);
\draw [ultra thick,red] (v1)--(v2)--(v4)--(v1);
\draw [ultra thick,red] (v4)--(v5)--(v7)--(v4);
\draw [ultra thick,red] (v7)--(v8)--(v10)--(v7);
\draw [ultra thick,red] (v10)--(v11)--(v1)--(v10);

\draw [ultra thick,blue] (v2)--(v3)--(v5)--(v2);
\draw [ultra thick,blue] (v5)--(v6)--(v8)--(v5);
\draw [ultra thick,blue] (v8)--(v9)--(v11)--(v8);
\draw [ultra thick,blue] (v11)--(v0)--(v2)--(v11);
\pgftransformxshift{7cm}
\pgftransformyshift{-0.5cm}
\foreach \x in {0,1,2} {
  \foreach \y in {0,1,2} {
    \node [vertex] (w\x\y) at (1.5*\x-1.5,1.5*\y-1.5) {};
  }
}
\draw (w00)--(w01)--(w02) [bend right = 30] to (w00);
\draw (w10)--(w11)--(w12) [bend right = 30] to (w10);
\draw (w20)--(w21)--(w22) [bend right = 30] to (w20);

\draw (w00)--(w10)--(w20) [bend left = 30] to (w00);
\draw (w01)--(w11)--(w21) [bend left = 30] to (w01);
\draw (w02)--(w12)--(w22) [bend left = 30] to (w02);
\node [vertex] (v12) at (3,0.75) {};
\node [vertex] (v02) at (0,3) {};
\node [vertex] (v01) at (-3,0.75) {};

\draw (v01)--(w00);
\draw (v01)--(w01);
\draw (v01)--(w02);
\draw (v01)--(w10);
\draw (v01)--(w11);
\draw (v01)--(w12);

\draw (v02)--(w00);
\draw (v02)--(w01);
\draw (v02)--(w02);
\draw (v02)--(w20);
\draw (v02)--(w21);
\draw (v02)--(w22);

\draw (v12)--(w10);
\draw (v12)--(w11);
\draw (v12)--(w12);
\draw (v12)--(w20);
\draw (v12)--(w21);
\draw (v12)--(w22);

\draw [ultra thick, red] (w00)--(w10)--(v01)--(w00);
\draw [ultra thick, red] (w01)--(w11)--(v01)--(w01);
\draw [ultra thick, red] (w02)--(w12)--(v01)--(w02);

\end{tikzpicture}
\end{center}
\caption{Two $12$-vertex exceptional graphs}
\label{fig:12vx}
\end{figure}
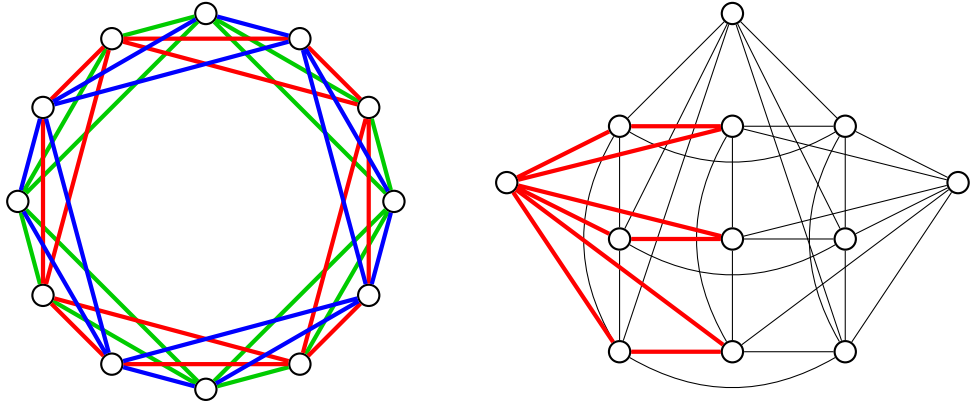

The first of the two $12$-vertex graphs is the circulant graph $\mathrm{Cay}(\mathbb{Z}_{12},\{\pm1,\pm2,\pm3\})$ which has a triangle decomposition containing the triangles
$\{\{a,a+1,a+3\} \mid a \in \mathbb{Z}_{12}\}$. This is depicted in the first diagram of \cref{fig:12vx} with a triangle decomposition.

The second of the two $12$-vertex graphs is obtained from the Cartesian product $K_3 \mathbin{\square} K_3$ by adjoining three vertices, $v_{01}$ which is adjacent to the six vertices $(x,y)$ where $x \in \{0,1\}$, then $v_{02}$ adjacent to the vertices with $x \in \{0,2\}$ and $v_{12}$ adjacent to the vertices with $x \in \{1,2\}$. This is depicted in the second diagram of \cref{fig:12vx} with a partial triangle decomposition.

The triangle decomposition contains three triangles from $K_3 \mathbin{\square} K_3$, namely the three triangles formed by vertices with the same first coordinate. Every other edge of $K_3 \mathbin{\square} K_3$ lies in a unique triangle with one of $\{v_{01},v_{02},v_{12}\}$.

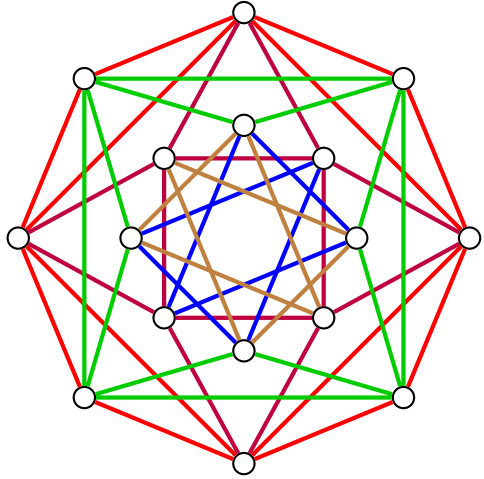
\begin{figure}
\begin{center}
\begin{tikzpicture}[scale=3]
\tikzstyle{vertex}=[fill=white, circle, draw=black, thick, inner sep=1mm]
\node [vertex] (v0) at (1.000000,0.000000) {};
\node [vertex] (v1) at (0.707107,0.707107) {};
\node [vertex] (v2) at (0.000000,1.000000) {};
\node [vertex] (v3) at (-0.707107,0.707107) {};
\node [vertex] (v4) at (-1.000000,0.000000) {};
\node [vertex] (v5) at (-0.707107,-0.707107) {};
\node [vertex] (v6) at (-0.000000,-1.000000) {};
\node [vertex] (v7) at (0.707107,-0.707107) {};
\node [vertex] (v8) at (0.500000,0.000000) {};
\node [vertex] (v9) at (0.353553,0.353553) {};
\node [vertex] (v10) at (0.000000,0.500000) {};
\node [vertex] (v11) at (-0.353553,0.353553) {};
\node [vertex] (v12) at (-0.500000,0.000000) {};
\node [vertex] (v13) at (-0.353553,-0.353553) {};
\node [vertex] (v14) at (-0.000000,-0.500000) {};
\node [vertex] (v15) at (0.353553,-0.353553) {};
\draw [ultra thick, red] (v0)--(v1)--(v2)--(v0);
\draw [ultra thick, red] (v2)--(v3)--(v4)--(v2);
\draw [ultra thick, red] (v4)--(v5)--(v6)--(v4);
\draw [ultra thick, red] (v6)--(v7)--(v0)--(v6);

\draw [ultra thick, purple] (v0)--(v9)--(v15)--(v0);
\draw [ultra thick, purple] (v6)--(v13)--(v15)--(v6);
\draw [ultra thick, purple] (v2)--(v9)--(v11)--(v2);
\draw [ultra thick, purple] (v4)--(v11)--(v13)--(v4);

\draw [ultra thick, green!80!black] (v1)--(v7)--(v8)--(v1);
\draw [ultra thick, green!80!black] (v1)--(v3)--(v10)--(v1);
\draw [ultra thick, green!80!black] (v3)--(v5)--(v12)--(v3);
\draw [ultra thick, green!80!black] (v5)--(v7)--(v14)--(v5);

\draw [ultra thick, blue] (v8)--(v10)--(v13)--(v8);
\draw [ultra thick, blue] (v9)--(v12)--(v14)--(v9);
\draw [ultra thick, brown] (v8)--(v11)--(v14)--(v8);
\draw [ultra thick, brown] (v10)--(v12)--(v15)--(v10);

\end{tikzpicture}
\label{fig:shrik}
\caption{Shrikhande graph with triangle decomposition}
\end{center}
\end{figure}

The first of the two $16$-vertex graphs is the well-known Shrikhande graph, which  has a triangle decomposition as illustrated in \cref{fig:shrik} (the colours in this figure have no meaning other than to identify the triangles). The second of the two graphs has a small automorphism group of order only $6$ and no particular structural features that would permit a simple description. 

Finally we move on to the main case, which is to consider which of the infinite number of linegraphs might occur as 2-distance graphs of cubic bipartite graphs of girth at least $6$. 

It is known that a connected linegraph $L(Y)$ of a graph $Y$ is regular if and only if $Y$ itself is connected and is a regular graph or is a \emph{bipartite $(a,b)$-semiregular} graph (this means that the vertices have degree either $a$ or $b$ according to which part of the bipartition they belong). The degree of $e = uv \in E(Y)$ in the linegraph $L(Y)$, is $d(u)+ d(v)-2$.  As $D_2(G)$ is $6$-regular, if $D_2(G)$ is a linegraph, then it is the linegraph of either a $4$-regular graphs or a bipartite $(a,b)$-semiregular graphs with $a+b=8$.

\begin{theorem}\label{tridecomp}
Let $Y$ be a connected $4$-regular graph. Then its linegraph $L(Y)$ has a triangle decomposition if and only if $Y = K_5$.
\end{theorem}

\begin{figure}
\begin{center}
\begin{tikzpicture}[scale=1.4]
\tikzstyle{vertex}=[fill=white, circle, ultra thick, draw=gray, inner sep=1mm]
\foreach \x in {0,1,...,9} {
\node [vertex] (v\x) at (36*\x:2cm) {};
}
\node (w13) at (0:2.3cm) {\footnotesize 13};
\node (w14) at (36*1:2.3cm) {\footnotesize 14};
\node (w04) at (36*2:2.3cm) {\footnotesize 04};
\node (w01) at (36*3:2.3cm) {\footnotesize 01};
\node (w12) at (36*4:2.3cm) {\footnotesize 12};
\node (w02) at (36*5:2.3cm) {\footnotesize 02};
\node (w03) at (36*6:2.3cm) {\footnotesize 03};
\node (w23) at (36*7:2.3cm) {\footnotesize 23};
\node (w24) at (36*8:2.3cm) {\footnotesize 24};
\node (w34) at (36*9:2.3cm) {\footnotesize 34};
\draw [ultra thick, red] (v0)--(v1)--(v3)--(v4)--(v0);
\draw [ultra thick, red] (v0)--(v3);
\draw [ultra thick, red] (v1)--(v4);

\draw [ultra thick, purple] (v4)--(v5)--(v7)--(v8)--(v4);
\draw [ultra thick, purple] (v4)--(v7);
\draw [ultra thick, purple] (v5)--(v8);

\draw [ultra thick, blue] (v2)--(v3)--(v5)--(v6)--(v2);
\draw [ultra thick, blue] (v2)--(v5);
\draw [ultra thick, blue] (v3)--(v6);

\draw [ultra thick, orange] (v6)--(v7)--(v9)--(v0)--(v6);
\draw [ultra thick, orange] (v6)--(v9);
\draw [ultra thick, orange] (v7)--(v0);

\draw [ultra thick, green!80!black] (v8)--(v9)--(v1)--(v2)--(v8);
\draw [ultra thick, green!80!black] (v8)--(v1);
\draw [ultra thick, green!80!black] (v9)--(v2);

\pgftransformxshift{6cm}

\foreach \x in {0,1,...,9} {
\node [vertex] (v\x) at (36*\x:2cm) {};
}
\node (w13) at (0:2.3cm) {\footnotesize 13};
\node (w14) at (36*1:2.3cm) {\footnotesize 14};
\node (w04) at (36*2:2.3cm) {\footnotesize 04};
\node (w01) at (36*3:2.3cm) {\footnotesize 01};
\node (w12) at (36*4:2.3cm) {\footnotesize 12};
\node (w02) at (36*5:2.3cm) {\footnotesize 02};
\node (w03) at (36*6:2.3cm) {\footnotesize 03};
\node (w23) at (36*7:2.3cm) {\footnotesize 23};
\node (w24) at (36*8:2.3cm) {\footnotesize 24};
\node (w34) at (36*9:2.3cm) {\footnotesize 34};

\draw [ultra thick, green!80!black] (v2)--(v8)--(v9)--(v2);
\draw [thick, gray] (v1)--(v4)--(v8)--(v1);
\draw [thick, gray] (v3)--(v4)--(v5)--(v3);
\draw [thick, gray] (v0)--(v4)--(v7)--(v0);
\draw [thick, gray] (v0)--(v3)--(v6)--(v0);
\draw [thick, gray] (v0)--(v1)--(v9)--(v0);
\draw [ultra thick, blue] (v2)--(v6)--(v5)--(v2);
\draw [thick,gray] (v1)--(v2)--(v3)--(v1);
\draw [ultra thick, orange] (v6)--(v7)--(v9)--(v6);
\draw [ultra thick, purple] (v5)--(v7)--(v8)--(v5);
\end{tikzpicture}
\end{center}
\caption{Edge-decompositions of $L(K_5)$ into $4$-cliques and into triangles}
\label{edgedecomp}
\end{figure}
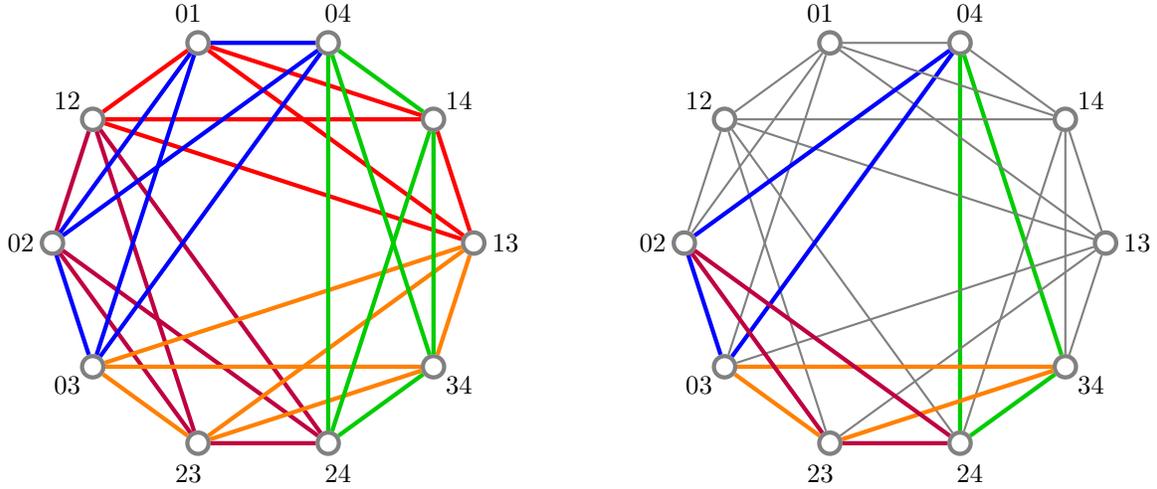

\begin{proof}
If $v \in V(Y)$ has neighbours $\{a,b,c,d\}$, then the four vertices $\{va, vb, vc, vd\}$ of $L(Y)$ form a $4$-clique. Each edge of $L(Y)$ is contained in exactly one of the $4$-cliques and so the edges of $L(Y)$ are partitioned into $4$-cliques---this is called a \emph{$K_4$-decomposition} of $L(Y)$. In what follows, we imagine the six edges of each $4$-clique as being \emph{coloured} with the name of the corresponding vertex of $Y$. For $Y=K_5$, this is illustrated in the first image of Figure~\ref{edgedecomp}.

For what follows, it helps to first consider the situation when $L(Y)$ \emph{does} have a triangle decomposition. The second image of Figure~\ref{edgedecomp} shows a triangle decomposition of $L(K_5)$. In this decomposition the triangles come in two types---there are \emph{monochromatic} triangles such as $\{04,24,34\}$ which lie entirely within a $4$-clique (these are highlighted in \cref{edgedecomp} with the corresponding colour), and \emph{colourful} triangles such as $\{01,12,02\}$ which contain edges from three different cliques (but are drawn in grey in \cref{edgedecomp}). Any monochromatic triangle has the form $\{va,vb,vc\}$ and any colourful triangle has the form $\{ab, bc, ca\}$, where $v$, $a$, $b$ and $c$ are vertices of $Y$. (The grey edges each lie in a unique grey triangle, so despite appearances, \cref{edgedecomp} completely determines the triangle decomposition.)

It is clear that each $4$-clique of $L(Y)$ contains at most one monochromatic triangle in a putative triangle decomposition and so there must be enough colourful triangles to cover the remaining edges. Our aim is to show that when $Y$ is not a complete graph, there are some edges in $L(Y)$ that cannot be incorporated into a suitable triangle.

We do this through an analysis of the neighbourhoods of the vertices in $Y$ identifying certain configurations in the neighbourhood of a single vertex that are incompatible with the existence of a triangle decomposition of $L(Y)$. For any vertex $v$ let $Y_1(v)$ denote the subgraph of $Y$ induced by the four vertices adjacent to $v$.

The following five cases enumerate various possibilities for $Y_1(v)$ that are incompatible with the existence of a triangle decomposition of $L(Y)$.

\begin{enumerate}[wide, label={\bf Case \arabic*:},ref={Case \arabic*}] 
\item  \label{noned} For some $v$, the graph $Y_1(v)$ has two vertex-disjoint non-edges.

Suppose that $a \not\sim b$ and $c \not\sim d$ in $Y_1(v)$. Then the edge $(va,vb)$ in $L(Y)$ cannot belong to a colourful triangle because the ``third vertex'' $ab$ is missing. So in any putative triangle decomposition of $L(Y)$, it must be the case that $(va,vb)$ lies in a monochromatic triangle of colour $v$. But the same argument implies that $(vc,vd)$ must also be contained in a monochromatic triangle of colour $v$, which is impossible.

\item  \label{K3isolated} For some $v$, the graph $Y_1(v)$ is isomorphic to $K_3 \cup K_1$.

Suppose that $\{a,b,c\}$ is a triangle in $Y_1(v)$ and $d$ is the isolated vertex. Then in any putative triangle decomposition of $L(Y)$, the edge $(va,vd)$ lies in a monochromatic triangle of colour $v$ and similarly for $(vb,vd)$ and $(vc,vd)$ which is again impossible.

\item \label{K13} For some $v$, the graph $Y_1(v)$ is isomorphic to $K_{1,3}$.

Suppose that $a$ is the vertex of degree $3$ in $Y_1(v)$ and that it is adjacent to $b$, $c$ and $d$. Then $b$ is adjacent to two more vertices $b'$ and $b''$ which are not equal to $c$ or $d$ and are not adjacent to $v$ or $a$. So $ab'$ and $vb''$ are two non-edges in $Y_1(b)$ and \ref{noned} shows that $L(Y)$ does not have a triangle decomposition.  

\item \label{K3pendant}For some $v$, the graph $Y_1(v)$ is isomorphic to $K_3$ with one pendant vertex.

Suppose that $\{a,b,c\}$ is a triangle in $Y_1(v)$ and $ad$ is the pendant edge. Then $d$ is adjacent two vertices $d'$ and $d''$ that are not equal to $b$ or $c$ and are not adjacent to $a$ or $v$. So $ad'$ and $vd''$ are two non-edges in $Y_1(d)$ and \ref{noned} shows that $L(Y)$ does not have a triangle decomposition.

\item \label{K4minuse} For some $v$, the graph $Y_1(v)$ is isomorphic to $K_4 - e$.

Suppose that every pair of vertices from $\{a,b,c,d\}$ is an edge of $Y_1(v)$ except for $cd$.  Then $c$ is adjacent to a vertex $c'$ that is not equal to $d$ and not adjacent to any of $a$, $b$ and $v$. Then $Y_1(c)$ is isomorphic to $K_3 \cup K_1$ and \ref{K3isolated} shows that $L(Y)$ does not have a triangle decomposition. 

\end{enumerate}

Now we combine all of these facts to complete the proof.  If there is a vertex $v$ for which $Y_1(v)$ has fewer than three edges, then \ref{noned} shows that $L(Y)$ has no triangle decomposition. If there is a vertex $v$ for which $Y_1(v)$ has exactly three edges then $Y_1(v)$ is either $P_3$, $K_3 \cup K_1$ or $K_{1,3}$ and \ref{noned}, \ref{K3isolated} or \ref{K13} (respectively) show that $L(Y)$ has no triangle decomposition. If there is a vertex $v$ for which $Y_1(v)$ has exactly four edges then $Y_1(v)$ is either $C_4$ or $K_3$ with a pendant vertex and \ref{noned} or \ref{K3pendant} (respectively) show that $L(Y)$ has no triangle decomposition.  If there is a vertex $v$ such that $Y_1(v)$ has exactly five edges, then $Y_1(v) = K_4-e$ and \ref{K4minuse} shows that $L(Y)$ does not have a triangle decomposition.

The only case remaining is when there is a vertex $v$ such that $Y_1(v)$ has six edges, in which case $Y = K_5$. Figure~\ref{edgedecomp} shows one triangle decomposition for $L(K_5)$ (the grey edges are in colourful triangles which are uniquely determined by a single edge).
\end{proof}

\begin{theorem}
Let $Y$ be a simple connected $(a,b)$-biregular graph with $a<b$ and $a+b=8$. Then its linegraph $L(Y)$ has a triangle decomposition if and only if $Y = K_{1,7}$.
\end{theorem}

\begin{proof}
Analogously to the situation in Theorem~\ref{tridecomp}, the edges of $L(Y)$ can be partitioned into cliques with some of size $a$ and some of size $b$. 

A colourful triangle in $L(Y)$ must have the form $\{ab, bc, ca\}$ and hence can only arise if $\{a,b,c\}$ are mutually adjacent in $Y$. Therefore when $Y$ is bipartite there are no colourful triangles in $L(Y)$ so a triangle decomposition has only monochromatic triangles. Therefore this can occur if and only if $K_a$ and $K_b$ separately have triangle decompositions. 

Neither $K_2$ (too few edges) nor $K_5$ (wrong number of edges modulo 3) have triangle decompositions and so $(a,b) \not= (2,6)$ and $(a,b) \not= (3,5)$. 
This only leaves $(a,b) = (1,7)$ in which case $Y = K_{1,7}$ and $L(Y)$ is a monochromatic $K_7$. Up to isomorphism, $K_7$ has a unique triangle decomposition into the seven blocks of the Fano plane.
\end{proof}

%% file: nonbipartite.tex
If $G$ is a non-bipartite cubic graph with no eigenvalues in $(-1,1)$ then its bipartite double $G \times K_2$ is either a \gm{} graph or one of the sporadic bipartite graphs listed in \cref{tab:maintable}.

If $G \times K_2$ is one of the sporadic graphs, then $G$ has at most $16$ vertices. As there are only $4681$ cubic graphs on $4$--$16$ vertices, direct computation is very straightforward and gives the results listed in \cref{tab:maintable}. Although there is only a minuscule chance that numerical approximation may cause eigenvalues very close to $\pm 1$ to be misclassified, we avoided even this small chance by using the function \verb+polsturm+ in Pari/GP \cite{PARI2}. This function uses exact arithmetic to count the number of real roots of an integer polynomial that lie in an interval given with exact end-points. 

First we consider the graphs whose bipartite double is the \gm{} graph.

\begin{theorem}\label{thm:KS-GM-bip-double}
If $G \times K_2$ is isomorphic to $\GM{k}$ then $k$ is even and $G$ is the \ks{} graph $\KS{k/2}$.
\end{theorem}

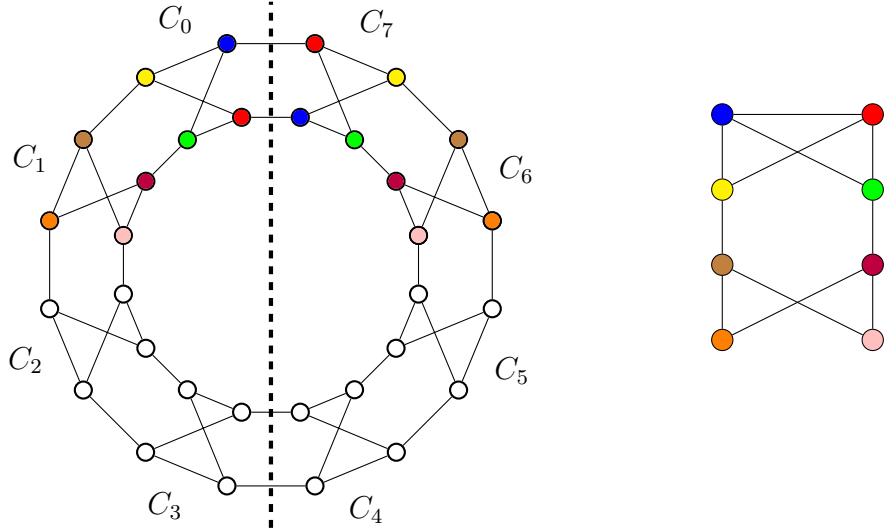
\begin{figure}
\begin{center}
\begin{tikzpicture}
\tikzstyle{vertex}=[fill=white, circle, draw=black, thick, inner sep=0.8mm]
\foreach \x in {0,1,...,16} {
 \node [vertex] (v\x) at (11.25 + 22.5*\x:3cm) {};
 \node [vertex] (w\x) at (11.25 + 22.5*\x:2cm) {};
}
\foreach \x in {0,1,...,7} {
  \node at (111.25+\x*45:3.5cm) {$C_\x$};
}
\draw (v0)--(w1)--(w0)--(v1)--(v0);
\draw (v2)--(w3)--(w2)--(v3)--(v2);
\draw (v4)--(w5)--(w4)--(v5)--(v4);
\draw (v6)--(w7)--(w6)--(v7)--(v6);
\draw (v8)--(w9)--(w8)--(v9)--(v8);
\draw (v10)--(w11)--(w10)--(v11)--(v10);
\draw (v12)--(w13)--(w12)--(v13)--(v12);
\draw (v14)--(w15)--(w14)--(v15)--(v14);
\draw (v1)--(v2);
\draw (v3)--(v4);
\draw (v5)--(v6);
\draw (v7)--(v8);
\draw (v9)--(v10);
\draw (v11)--(v12);
\draw (v13)--(v14);
\draw (v15)--(v0);

\draw (w1)--(w2);
\draw (w3)--(w4);
\draw (w5)--(w6);
\draw (w7)--(w8);
\draw (w9)--(w10);
\draw (w11)--(w12);
\draw (w13)--(w14);
\draw (w15)--(w0);

\draw [dashed, ultra thick] (0,-3.5)--(0,3.5);

\fill [red] (v3) circle (1mm);
\fill [red] (w4) circle (1mm);

\fill [blue] (w3) circle (1mm);
\fill [blue] (v4) circle (1mm);

\fill [yellow] (v2) circle (1mm);
\fill [yellow] (v5) circle (1mm);

\fill [green] (w2) circle (1mm);
\fill [green] (w5) circle (1mm);

\fill [purple] (w1) circle (1mm);
\fill [purple] (w6) circle (1mm);

\fill [brown] (v1) circle (1mm);
\fill [brown] (v6) circle (1mm);

\fill [orange] (v0) circle (1mm);
\fill [orange] (v7) circle (1mm);

\fill [pink] (w0) circle (1mm);
\fill [pink] (w7) circle (1mm);

\pgftransformxshift{7cm}

\node [fill=blue,draw=black,circle,inner sep=1mm] (b) at (-1,2) {};
\node [fill=red,draw=black,circle,inner sep=1mm] (r) at (1,2) {};
\node [fill=yellow,draw=black,circle,inner sep=1mm] (y) at (-1,1) {};
\node [fill=green,draw=black,circle,inner sep=1mm] (g) at (1,1) {};
\node [fill=purple,draw=black,circle,inner sep=1mm] (p) at (1,0) {};
\node [fill=brown,draw=black,circle,inner sep=1mm] (br) at (-1,0) {};
\node [fill=orange,draw=black,circle,inner sep=1mm] (or) at (-1,-1) {};
\node [fill=pink,draw=black,circle,inner sep=1mm] (pi) at (1,-1) {};
\draw (b)--(r);
\draw (b)--(y);
\draw (b)--(g);
\draw (y)--(br);
\draw (y)--(r);
\draw (g)--(p);
\draw (g)--(r);
\draw (br)--(or);
\draw (br)--(pi);
\draw (p)--(or);
\draw (p)--(pi);
\end{tikzpicture}
\end{center}
\caption{$\GM{8} = \KS{4} \times K_2$}
\label{fig:sigma}
\end{figure}

\begin{proof}
Imrich and Pisanski \cite{MR2419215} show that a bipartite graph is a \emph{Kronecker cover} (i.e., of the form $G \times K_2$) if and only if its automorphism group has a fixed-point-free involution $\sigma$ that exchanges the colour classes and is such that $v \not \sim v^\sigma$ for all $v$. In this case, the orbits of the involution are the fibres of the covering map and so the quotient graph can easily be recovered. We will prove this theorem by showing that, without loss of generality, any suitable automorphism $\sigma$ has the \ks{} graph on half the number of vertices as quotient.

The automorphism group of $\GM{k}$ is the semidirect product $\mathbb{Z}_2^k \rtimes D_k$ where $D_k$ is the dihedral group of order $2k$. The elementary abelian subgroup $\mathbb{Z}_2^k$ is generated by the $k$ elements that swap the two edges in the matching that connects adjacent $4$-cycles while fixing all other vertices. The dihedral group acts on the $k$ induced $4$-cycles by rotation and/or reflection. 

In $\GM{k}$ a colour-swapping involution that does not fix any edges cannot fix any of the $4$-cycles and therefore $k$ must be even. \cref{fig:sigma} illustrates the partially-completed construction of such an involution $\sigma$ for $k=8$, using colours to indicate which pairs of vertices are exchanged by the involution.  We may assume without loss of generality that the reflection from the dihedral group has the vertical axis shown in \cref{fig:sigma}, in which case $\sigma$ must exchange the two blue vertices and the two red vertices. Now $\sigma$ must exchange the other two vertices in $C_0$ with the other two vertices in $C_7$. Because these pairs of vertices are completely joined to the other vertices in the $4$-cycle, we can assume without of generality that the green vertices are paired and the yellow vertices are paired. Now the two vertices adjacent to the two yellow vertices must be exchanged by $\sigma$; these are coloured brown in the diagram, and similarly the two purple vertices each adjacent to one of the two green vertices must be exchanged. Continuing in this fashion, the involution $\sigma$ is constructed by a series of choices where each choice is either forced, or can be made without loss of generality. The second diagram shows the partially-constructed quotient graph that is emerging from this process, which can only be completed to $\KS{4}$. 
\end{proof}

%% file: maximalgapset.tex
\subsection{A maximal gap set}

Koll\'ar and Sarnak \cite{KolSar2021} prove that $(-1,1)$ is a maximal gap \emph{interval} but raise the question of whether it is a maximal gap \emph{set}. In other words, are there any \emph{additional} open intervals that are avoided by infinitely many graphs chosen from those that already avoid $(-1,1)$?

From our main result, any such infinite set of graphs must contain infinitely many Guo-Mohar or Koll\'ar-Sarnak graphs.  The next lemma shows that in terms of eigenvalues we may assume that such a set contains only Guo-Mohar graphs.

\begin{lemma}\label{sameeigs}\
The characteristic polynomials of the Guo-Mohar graphs and the Koll\'ar-Sarnak graphs are related by the equation
\[
\varphi (\GM{2k},x)  = 
\left( \frac{x+3} {x-3} \right)
\left( \frac{x-1} {x+1} \right)^3
\varphi( \KS{k}, x )^2.
\]
\end{lemma}

\begin{proof}
The characteristic polynomials of the Guo-Mohar and Koll\'ar-Sarnak graphs satisfy recurrence relations that can be solved explicitly and then used to verify this expression. The full details are somewhat tedious and not very illuminating and therefore we omit them here. \qedhere
\end{proof}

The Guo-Mohar graph $\GM{k}$ on $4k$ vertices has eigenvalues $\pm 1$ each with multiplicity $k$ and the values
\[
\tau_j = \pm 
\sqrt{
5 + 4 \cos \left( 2 \pi j / k \right)
}
\]
for $0 \leqslant j \leqslant k-1$. (This expression may increase the stated multiplicity of $\pm 1$.)  

The eigenvalues of $\GM{k}$ are the images of the rational points $j/k$ for $0 \leqslant j \leqslant k-1$ for under the map 
\[
x \mapsto \sqrt{5 + 4\cos (\pi x)}.
\]
If we consider this map as a function with domain $[0,1]$ then it has range $[1,3]$ and is continuous and decreasing. Therefore any open interval $\mathcal{I} \subset [1,3]$ has some preimage $\mathcal{I'} \subset [0,1]$ which is also an open interval. If $1/K$ is less than the length of this preimage then $\mathcal{I'}$ contains some multiple of $1/k$ for every $k > K$. We conclude that $\mathcal{I}$ contains eigenvalues of $\GM{k}$ for all $k > K$.  A symmetric argument holds for open intervals contained in $[-3,-1]$.

\begin{theorem}
The interval $(-1,1)$ is a maximal gap set for cubic graphs.
\end{theorem}

\begin{proof}
By \cref{sameeigs} every Kolla\'r-Sarnak graph has precisely the same eigenvalues as a Guo-Mohar graph (with the sole exception of $-3$). Therefore in any infinite set of cubic graphs avoiding $(-1,1)$ we may replace each Koll\'ar-Sarnak graph with a Guo-Mohar graph without altering the overall set of eigenvalues. But from the discussion in the preceding paragraph, every open interval contained in $[-3,-1]$ or $[1,3]$ contains eigenvalues from all sufficiently large Guo-Mohar graphs and hence cannot be a subset of a gap set. \qedhere
\end{proof}

%% file: conclusion.tex
\section{Open problems and future work}

Lower and upper bounds on the maximum value of the HL-index among all graphs with given average degree are given in \cite{Moh2015}. In the same paper, it is shown that a positive fraction of the eigenvalues of a subcubic graph lie in the interval $[-\sqrt{2}, \sqrt{2}]$.
It is still an open problem to show that the median eigenvalues of any cubic (or subcubic) graph, apart from the Heawood graph, are in the interval $[-1,1]$. 

An open sub-problem posed in \cite{Moh2013} is to show that the median eigenvalues of subcubic planar graphs lie in $[-1,1]$; this problem was also collected in the recent collection of open problems in spectral graph theory, see \cite[Section 19.2]{LiuNin2023}. It is true for subcubic, planar and  $K_4$-free graphs, as shown in \cite{WanZha2024}. 

The proofs in this paper rely heavily on the graph being cubic, and so do not apply to non-regular subcubic graphs. However we suspect that a similar result is true for non-regular subcubic graphs. In particular, we know two infinite families of non-regular subcubic graphs that arise from the Koll\'ar-Sarnak graphs. The first family is obtained from $\KS{k}$ by deleting $w_0$, and second is obtained from $\KS{k}$ by deleting $\{w_0,b_{k-1}\}$ (see \cref{ksfig} for the vertex-naming convention). In addition we know four ``sporadic'' graphs on $8$, $10$, $14$ and $18$ vertices but we cannot currently rule out the existence of larger ones.

\section*{Acknowledgements} 

We would like to thank Brendan McKay who first pointed out the connection to positive semi-definite matrices and whose computations provided independent verification that the list of examples on up to $32$ vertices is complete. 

K.~Guo gratefully acknowledge the support of the Cheryl E.~Praeger Visiting Research Fellowship, which facilitated the initiation of this research during a visit to the University of Western Australia.